\documentclass{amsart}
\usepackage{amssymb,latexsym}
\usepackage{amsmath}
\usepackage{graphicx}
\usepackage{amscd}
\usepackage{color}
\usepackage{enumerate}
\newenvironment{enumeratei}{\begin{enumerate}[\upshape (i)]}{\end{enumerate}}
\numberwithin{equation}{section}
\theoremstyle{plain}
 \newtheorem{theorem}{Theorem}[section]

\theoremstyle{definition}

\newcommand \datum {\hfill  December 17, 2017}
\newcommand \con {\textup{con}}
\newcommand \Con  {\textup{Con}}
\newcommand\ideal[1]{\mathord\downarrow #1}
\newcommand\filter[1]{\mathord\uparrow #1}
\newcommand \blokk [2] {#1/#2}

\newcommand \tbf [1] {\textbf{#1}} 
\newcommand \set[1] {\{#1\}}
\newcommand \tuple [1] {\langle #1\rangle}
\newcommand \pair [2] {\tuple{#1,#2}}

\newcommand \url [1] {\texttt{#1}} 

\newcommand \red [1] {\color{red}#1\color{black}}

\begin{document}
\title
[Finite lattices with many congruences]
{A note on finite lattices with many congruences}

\author[G.\ Cz\'edli]{{G\'abor Cz\'edli}}
\email{czedli@math.u-szeged.hu}
\urladdr{http://www.math.u-szeged.hu/~czedli/}
\address{Bolyai Institute, University of Szeged, Hungary 6720}

\begin{abstract} By a twenty year old result of Ralph Freese, an $n$-element lattice  $L$ has at most $2^{n-1}$ congruences. We prove that if $L$ has less than $2^{n-1}$ congruences, then it has at most $2^{n-2}$ congruences. Also, we describe the $n$-element lattices with exactly $2^{n-2}$ congruences. 
\end{abstract}


\maketitle
\section{Introduction and motivation}
It follows from  Lagrange's Theorem that the size $|S|$ of an arbitrary subgroup $S$ of a finite group $G$ is either $|G|$, or it is at most the half of the maximum possible value, $|G|/2$. Furthermore, if the size of $S$ is the half of its maximum possible value, then $S$ has some special property since it is normal. Our goal is to prove something similar on the size of the congruence lattice $\Con(L)$ of an $n$-element lattice $L$. 

For a \emph{finite} lattice $L$, the relation between $|L|$ and $|\Con(L)|$ has been studied in some earlier papers, including Freese~\cite{freesecomp}, 
Gr\"atzer and Knapp~\cite{grknapp4}, 
Gr\"atzer, Lakser, and Schmidt~\cite{GGHLSchT}, 
Gr\"atzer, Rival, and Zaguia~\cite{GRivalZ}. In particular,  part~\eqref{thmmaina} of Theorem~\ref{thmmain} below is due to Freese~\cite{freesecomp}; note that we are going present a new proof of part~\eqref{thmmaina}.
Although Cz\'edli and Mure\c san~\cite{czgmuresan} and  Mure\c san~\cite{muresan2017arXiv} deal only with infinite lattices, they are also among the papers  motivating the present one.

\section{Our result and its proof}
Mostly, we  follow the terminology and notation of Gr\"atzer~\cite{ggfoundbook}. In particular, the \emph{glued sum}  $L_0\dot+ L_1$ 
of finite lattices $L_0$  and 
$L_1$ is their Hall--Dilworth gluing along $L_0\cap L_1=\set{1_{L_0}}=\set{0_{L_1}}$; see, for example, Gr\"atzer~\cite[Section IV.2]{ggfoundbook}. Note that $\dot+$ is an associative operation. Our result is the following.

\begin{theorem}\label{thmmain} If $L$ is a finite lattice of size $n=|L|$, then the following hold.
\begin{enumeratei}
\item\label{thmmaina} $L$ has at most $2^{n-1}$ many congruences. Furthermore, $|\Con(L)|=2^{n-1}$  if and only if $L$ is a chain. 
\item\label{thmmainb} If $L$ has \emph{less than} $2^{n-1}$  congruences, then it has at most $2^{n-1}/2 = 2^{n-2}$  congruences.
\item\label{thmmainc} $|\Con(L)|=2^{n-2}$ if and only if $L$ is  of the form  $C_1\dot+B_2\dot+ C_2$ such that $C_1$ and $C_2$ are chains and $B_2$ is the four-element Boolean lattice.
\end{enumeratei}
\end{theorem}

\begin{figure}[ht] 
\centerline
{\includegraphics[scale=1.0]{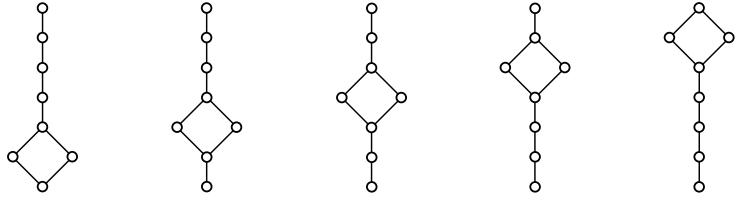}}
\caption{The full list of 8-element lattices with exactly $64=2^{8-2}$ many congruences}
\label{figone}
\end{figure}

For $n=8$, part \eqref{thmmainc} of this theorem is illustrated in Figure~\ref{figone}. Note that part \eqref{thmmaina} of the theorem is due to Freese~\cite[page 3458]{freesecomp}; however, our approach to Theorem~\ref{thmmain} includes a new proof of part \eqref{thmmaina}. 

\begin{proof}[Proof of Theorem~\ref{thmmain}]
We prove the theorem by induction on $n=|L|$. Since the case $n=1$ is clear, assume as an induction hypothesis that $n>1$ is a natural number and all the three parts of the theorem hold for every lattice with size less than $n$. Let $L$ be a lattice with $|L|=n$. For $\pair a b\in L^2$, the least congruence collapsing $a$ and $b$ will be denoted by $\con(a, b)$. 
A \emph{prime interval} or an \emph{edge} of $L$ is an interval $[a,b]$ with $a\prec b$. For later reference, note that
\begin{equation}
\parbox{7.5cm}{$\Con(L)$ has an atom, and every of its atoms is of the form $\con(a,b)$ for some
prime interval $[a,b]$;}
\label{eqtxtmTmsdsvb}
\end{equation}
this follows from the finiteness of $\Con(L)$ and from the fact that every congruence on $L$ is the join of congruences generated by \emph{covering pairs} of elements; see also Gr\"atzer~\cite[page 39]{ggCFL2} for this folkloric fact.

Based on \eqref{eqtxtmTmsdsvb}, pick a prime interval $[a,b]$ of $L$ such that $\Theta=\con(a,b)$ is  an atom in $\Con(L)$. Consider the map $f\colon \Con(L)\to \Con(L)$ defined by $\Psi\mapsto \Theta\vee \Psi$. We claim that, with respect to $f$,  
\begin{equation}
\text{every element of $f(\Con(L))$ has at most two  preimages.}
\label{eqtxttmostttT}
\end{equation}
Suppose to the contrary that there are pairwise distinct $\Psi_1,\Psi_2,\Psi_3\in\Con(L)$ with the same $f$-image.
Since the $\Theta \wedge \Psi_i$ belong to the two-element principal ideal $\ideal \Theta:=\set{\Gamma\in \Con(L):\Gamma\leq \Theta}$ of $\Con(L)$, at least two of these meets coincide. So we can assume that $\Theta\wedge \Psi_1=
\Theta\wedge \Psi_2$ and, of course, we have that $\Theta\vee \Psi_1=f(\Psi_1)=f(\Psi_2) = \Theta\vee \Psi_2$. This means that both $\Psi_1$ and $\Psi_2$ are relative complements of $\Theta$ in the interval $[\Theta\wedge \Psi_1,\Theta\vee \Psi_1]$. 
According to a classical result of Funayama and Nakayama~\cite{funayamanakayama}, $\Con(L)$ is distributive. Since relative complements in distributive lattices are well-known to be unique, see, for example, Gr\"atzer~\cite[Corollary 103]{ggfoundbook}, it follows that $\Psi_1=\Psi_2$. This is a contradiction proving \eqref{eqtxttmostttT}.

Clearly, $f$ is a retraction map onto the filter $\filter \Theta$. It follows from \eqref{eqtxttmostttT} that $|\filter\Theta|\geq |\Con(L)| / 2$. Also, by the well-known Correspondence Theorem, see Burris and Sankappanawar~\cite[Theorem 6.20]{bursan}, or see  Theorem 5.4 (under the name Second Isomorphism Theorem) in Nation~\cite{nationbook}, 
$|\filter\Theta| = |\Con(L/\Theta)|$ holds. Hence, it follows that
\begin{equation}
|\Con(L/\Theta)|\geq \frac 12\cdot |\Con(L)|.
\label{eqdlshbzrGt}
\end{equation}
Since $\Theta$ collapses at least one pair of distinct elements, $\pair a b$, we have that $|L/\Theta|\leq n-1$. 
Thus, it follows from part \eqref{thmmaina} of the induction hypothesis  that $|\Con(L/\Theta)|\leq 2^{(n-1)-1}=2^{n-2}$. Combining this inequality with  \eqref{eqdlshbzrGt}, we obtain that $|\Con(L)|\leq 2\cdot |\Con(L/\Theta)| \leq 2^{n-1}$. This shows the first half of part \eqref{thmmaina}.    

If $L$ is a chain, then $\Con(L)$ is known to be the $2^{n-1}$-element boolean lattice; see, for example,
Gr\"atzer~\cite[Corollaries 3.11 and 3.12]{ggCFL2}. Hence, we have that $|\Con(L)|=2^{n-1}$ if $L$ is a chain. Conversely, 
assume the validity of $|\Con(L)|=2^{n-1}$, 
and let $k=|L/\Theta|$. By the induction hypothesis, $|\Con(L/\Theta)|\leq 2^{k-1}$. On the other hand, $|\Con(L/\Theta)|\geq |\Con(L)|/2=2^{n-2}$ holds by \eqref{eqdlshbzrGt}.
These two inequalities and $k<n$ yield that $k=n-1$ and also that $|\Con(L/\Theta)|=2^{n-2}=2^{k-1}$. Hence, the induction hypothesis implies that $L/\Theta$ is a chain. For the sake of contradiction, suppose that $L$ is not a chain, and pick a pair $\pair u v$ of incomparable elements of $L$. The $\Theta$-blocks $\blokk u\Theta$ and $\blokk v\Theta$ are comparable elements of the chain $L/\Theta$, whence we can assume that $\blokk u\Theta\leq\blokk v\Theta$. It follows that $\blokk u\Theta =\blokk u\Theta \wedge \blokk v\Theta = \blokk {(u\wedge v)}\Theta$ and, by duality, $\blokk v\Theta  = \blokk {(u\vee v)}\Theta$. Thus, since $u$, $v$, $u\wedge v$ and $u\vee v$ are pairwise distinct elements of $L$ and $\Theta$ collapses both of the pairs $\pair {u\wedge v}u$ and $\pair{v}{u\vee v}$, we have that $k=|L/\Theta|<n-2$, which is a contradiction. This proves part \eqref{thmmaina} of the theorem.

As usual, for a lattice $K$, let $J(K)$ and $M(K)$ denote the set of nonzero \emph{join-irreducible elements} and the set of \emph{meet-irreducible elements} distinct from $1$, respectively. By a \emph{narrows} we will mean a prime interval $[a,b]$ such that $a\in M(L)$ and $b\in J(L)$. Using  Gr\"atzer~\cite{ggtechnolemma}, it follows in a straightforward way that 
\begin{equation}
\parbox{6cm}{if $[a,b]$ is a narrows, then $\set{a,b}$ is the only non-singleton block of $\con(a,b)$.}
\label{eqtxtnrrsws}
\end{equation}
Now, in order to prove part \eqref{thmmainb} of the theorem, assume that $|\Con(L)|<2^{n-1}$. By \eqref{thmmain}, we can pick a prime interval $[a,b]$ such that $\Theta:=\con(a,b)$ is an atom in $\Con(L)$. There are two cases to consider depending on whether $[a,b]$ is a narrows or not; for later reference, some parts of the arguments for these two cases will be summarized in  \eqref{eqtxtlsRfhrZnB} and \eqref{eqtxtdzNbvQxRtk} redundantly.
First, we deal with the case where  $[a,b]$ is a narrows. We claim that 
\begin{equation}
\parbox{8.7cm}{if $|\Con(L)|<2^{n-1}$, $[a,b]$ is a narrows, and $\Theta=\con(a,b)$ is an atom in $\Con(L)$, then $L/\Theta$ is not a chain.}
\label{eqtxtlsRfhrZnB}
\end{equation}
By \eqref{eqtxtnrrsws}, $|L/\Theta|=n-1$. By the already proved part \eqref{thmmaina}, $L$ is not a chain, whence there are $u,v\in L$ such that $u\parallel v$. We claim that $\blokk u\Theta$ and $\blokk v\Theta$ are incomparable elements of $L/\Theta$. Suppose the contrary. Since $u$ and $v$ play a symmetric role, we can assume that $\blokk u\Theta\vee \blokk v\Theta=\blokk v\Theta$, i.e., $\blokk{(u\vee v)}\Theta=\blokk v\Theta$.  But $u\vee v\neq v$ since $u\parallel v$, whereby \eqref{eqtxtnrrsws} gives that $\set{v,u\vee v}=\set{a,b}$. Since $a<b$, this means that $v=a$ and $u\vee v=b$. Thus, $u\vee v\in J(L)$ since $[a,b]$ is a narrows. The membership $u\vee v\in J(L)$ gives that $u\vee v\in\set{u,v}$, contradicting $u\parallel v$. This shows that $\blokk u\Theta\parallel\blokk v\Theta$, whence $L/\Theta$ is not a chain. We have shown the validity of \eqref{eqtxtlsRfhrZnB}.
Using part \eqref{thmmaina} and $|L/\Theta|=n-1$, it follows that $|\Con(L/\Theta)|<2^{(n-1)-1}$. By the induction hypothesis, we can apply \eqref{thmmainb}  to $L/\Theta$ to conclude that $|\Con(L/\Theta)|\leq 2^{(n-1)-2}$. This inequality and \eqref{eqdlshbzrGt} yield that $|\Con(L)|\leq 2\cdot |\Con(L/\Theta)|\leq 2^{n-2}$, as required.  

Second, assume that $[a,b]$ is not a narrows. Our immediate plan is to show that
\begin{equation}
\parbox{7.2cm}{if a prime interval $[a,b]$ of $L$ is not a narrows and
$\Theta=\con(a,b)$, then $|L/\Theta|\leq n-2$.}
\label{eqtxtdzNbvQxRtk}
\end{equation}
By duality, we can assume that $a$ is meet-reducible. Hence, we can pick an element $c\in L$ such that $a\prec c$ and $c\neq b$. 
Clearly, $c\neq b\vee c$ and $\Theta=\con(a,b)$ collapses both $\pair a b$ and $\pair c{b\vee c}$, which are distinct pairs. 
Thus, we obtain that $|L/\Theta|\leq n-2$, proving  \eqref{eqtxtdzNbvQxRtk}. Hence, $\Con(L/\Theta)\leq 2^{n-3}$ by part \eqref{thmmaina} of the induction hypothesis. Combining this inequality with \eqref{eqdlshbzrGt}, we obtain the required inequality $\Con(L)\leq 2^{n-2}$. This completes the induction step for part \eqref{thmmainb}.

\begin{figure}[ht] 
\centerline
{\includegraphics[scale=1.0]{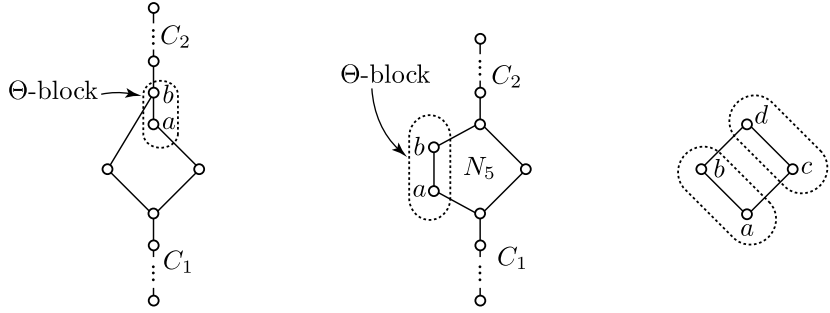}}
\caption{Illustrations for the proof}
\label{figtwo}
\end{figure}

Next, in order to perform the induction step for  part \eqref{thmmainc}, we assume that $|\Con(L)|=2^{n-2}$. Again, there are two cases to consider. First, we assume that there exists a narrows $[a,b]$ in $L$ such that $\Theta:=\con(a,b)$ is an atom in $\Con(L)$.  Then $|L/\Theta|=n-1$ by \eqref{eqtxtnrrsws} and $L/\Theta$ is not a chain by \eqref{eqtxtlsRfhrZnB}. By the  induction hypothesis, parts \eqref{thmmaina} and \eqref{thmmainb} hold for $L/\Theta$, whereby we have that $|\Con(L/\Theta)|\leq 2^{(n-1)-2}=2^{n-3}$. On the other hand, it follows from \eqref{eqdlshbzrGt} that $|\Con(L/\Theta)|\geq |\Con(L)|/2= 2^{n-3}$. Hence, 
$|\Con(L/\Theta)| =  2^{n-3}= 2^{|L/\Theta|-2}$. By the induction hypothesis,  $L/\Theta$ is of the form $C_1\dot+B_2\dot+ C_2$. We know from \eqref{eqtxtnrrsws} that $\set{a,b}=[a,b]$ is the unique non-singleton $\Theta$-block. If this $\Theta$-block is outside $B_2$, then $L$ is obviously of the required form. If the  $\Theta$-block $\set{a,b}$ is in $C_2\cap B_2$, then $L$ is of the required form simply because the situation on the left of Figure~\ref{figtwo} would contradict the fact that $[a,b]$ is a narrows. A dual treatment applies for the case $\set{a,b}\in C_1\cap B_2$. 
If the  $\Theta$-block $\set{a,b}$ is in $B_2\setminus{(C_1\cup C_2})$, then $L$ is of the form $C_1\dot+ N_5 \dot+ C_2$, where $N_5$ is the ``pentagon''; see the middle part of Figure~\ref{figtwo}. For an arbitrary  bounded lattice $K$ and the two-element chain $\mathbf 2$, it is straightforward to see that 
\begin{equation}
\Con(K\dot+\mathbf 2)\cong \Con(\mathbf 2\dot+K)\cong \Con(K)\times \mathbf 2.
\label{eqczNsPznGrHV}
\end{equation}
A trivial induction based on \eqref{eqczNsPznGrHV} yields that $|\Con(C_1\dot+ N_5 \dot+ C_2)|$ is divisible by $5=|\Con(N_5)|$. But $5$ does not divide $|\Con(L)|=2^{n-2}$, ruling out the case that the $\Theta$-block $\set{a,b}$ is in $B_2\setminus{(C_1\cup C_2})$. Hence, $L$ is of the required form.

Second, we assume that no narrows in $L$ generates an atom of $\Con(L)$. 
By \eqref{eqtxtmTmsdsvb}, we can pick a prime interval $[a,b]$ such that $\Theta:=\con(a,b)$ is an atom of $\Con(L)$. Since $[a,b]$ is not a narrows, 
\eqref{eqtxtdzNbvQxRtk} gives that $|L/\Theta|\leq n-2$. We claim that we have equality here, that is, $|L/\Theta|= n-2$. Suppose to the contrary that $|L/\Theta|\leq n-3$. Then part \eqref{thmmaina} and \eqref{eqdlshbzrGt} yield that 
\[2^{n-2}=|\Con(L)|\leq 2\cdot |\Con(L/\Theta)|\leq 2\cdot 2^{(n-3)-1}=2^{n-3},
\]
which is a contradiction. Hence,  $|L/\Theta|= n-2$. Thus, we obtain from by part \eqref{thmmaina} that $|\Con(L/\Theta)|\leq 2^{n-3}$.  On the other hand, 
 \eqref{eqdlshbzrGt} yields that  $|\Con(L/\Theta)|\geq |\Con(L)|/2=2^{n-3}$, whence  $|\Con(L/\Theta)|=2^{n-3}=2^{|L/\Theta|-1}$, and it follows by part \eqref{thmmaina} that $L/\Theta$ is a chain. Now, we have to look at the prime interval $[a,b]$ closely. It is not a narrows, whereby duality allows us to assume that $b$ is not the only cover of $a$. So we can pick an element $c\in L\setminus\set{b}$ such that $a\prec c$, and let $d:=b\vee c$; see on the right of Figure~\ref{figtwo}. Since $\pair c d=\pair{c\vee a}{c\vee b} \in \con(a,b)= \Theta$,  any two elements of $[c,d]$ is collapsed by $\Theta$. Using $\pair a b\in \Theta$, $\pair c d\in \Theta$, and $|L/\Theta|=n-2=|L|-2$, it follows that there is no ``internal element'' in the interval $[c,d]$, that is, $c\prec d$.  Furthermore, $[a,b]=\set{a,b}$ and $[c,d]=\set{c,d}$ are the only non-singleton blocks of $\Theta$. 
In order to show that $b\prec d$, suppose to the contrary that $b<e<d$ holds for some $e\in L$. Since $d=b\vee c\leq e\vee c \leq d$, we have that $e\vee c=d$, implying $c\nleq e$. Hence, $c\wedge e<e$. Since 
$\pair {c\wedge e}e = \pair {c\wedge e}{d\wedge e}\in\Theta$, the $\Theta$-block of $e$ is not a singleton. This contradicts the fact that  $\set{a,b}$ and $\set{c,d}$ are the only non-singleton $\Theta$-blocks, whereby we conclude that $b\prec d$. The covering relations established so far show that $S:=\set{a=b\wedge c, b, c, d=b\vee c}$ is a covering square in $L$. We know that both non-singleton $\Theta$-blocks are subsets of $S$ and $L/\Theta$ is a chain. Consequently, $L\setminus S$ is also a chain.

Hence, to complete the analysis of the second case when $[a,b]$ is not a narrows, it suffices to show that for every $x\in L\setminus S$, we have that either $x\leq a$, or  $x\geq d$. So, assume that $x\in L\setminus S$. Since $L/\Theta$ is a chain, $\set{a,b}$ and $\set x$ are comparable in $L/\Theta$. If $\set x < \set{a,b}$, then 
$\set x\vee \set{a,b}=\set{a,b}$ gives that $x\vee a\in \set{a,b}$. If $x\vee a$ happens to equal $b$, then $x\nleq a$ leads to $x\wedge a<x$ and 
$\pair{x\wedge a}x =\pair{x\wedge a}{x\wedge b}\in\Theta$, contradicting the fact the $\set{a,b}$ and $\set{c,d}$ are the only non-singleton $\Theta$-blocks. So if  $\set x < \set{a,b}$, then $x\vee a=a$ and $x<a$, as required. Thus, we can assume that $\set x > \set{a,b}$. If $\set x > \set{c,d}$, then the dual of the easy argument just completed shows that $x\geq d$. So, we are left with the case $\set{a,b}<\set x < \set{c,d}$. Then the equalities $\set{a,b}\vee \set x=\set x$ and $\set x=\set x\wedge \set{c,d}$ give that $b\vee x=x=d\wedge x$, that is, $b\leq x\leq d$. But $x\notin S$, so $b<x<d$, contradicting $b\prec d$. This completes the second case of the induction step for part \eqref{thmmainc} and the proof of Theorem~\ref{thmmain}. 
\end{proof}

\end{document}